\documentclass[11pt,a4paper]{article}

\usepackage{amsmath}
\usepackage{amsthm}
\usepackage{amssymb}
\usepackage{amscd}
\usepackage[all]{xy}

\title{A Note on the Effective Non-vanishing Conjecture}
\author{Qihong Xie}
\date{}
\theoremstyle{plain}
\newtheorem{prop}{Proposition}
\newtheorem{lem}[prop]{Lemma}
\newtheorem{thm}[prop]{Theorem}
\newtheorem{cor}[prop]{Corollary}
\newtheorem{conj}[prop]{Conjecture}

\theoremstyle{definition}

\newtheorem*{ack}{Acknowledgments}

\theoremstyle{remark}
\newtheorem{rem}[prop]{Remark}

\newcommand{\Q}{\mathbb Q}
\newcommand{\R}{\mathbb R}
\newcommand{\C}{\mathbb C}

\newcommand{\OO}{\mathcal O}

\newcommand{\GG}{\mathcal G}

\newcommand{\FF}{\mathcal F}

\newcommand{\al}{\alpha}
\newcommand{\Pic}{\mathop{\rm Pic}\nolimits}
\newcommand{\Alb}{\mathop{\rm Alb}\nolimits}
\newcommand{\rk}{\mathop{\rm rank}\nolimits}

\begin{document}

\maketitle

\begin{abstract}
We give a reduction of the irregular case for the effective 
non-vanishing conjecture by virtue of the Fourier-Mukai transform. 
As a consequence, we reprove that the effective non-vanishing 
conjecture holds on algebraic surfaces.
\end{abstract}

In this note we consider the following so-called 
effective non-vanishing conjecture, which has been 
put forward by Ambro and Kawamata \cite{am, ka00}.

\begin{conj}[$EN_n$]\label{1}
Let $X$ be a proper normal variety of dimension $n$, 
$B$ an effective $\R$-divisor on $X$ such that 
the pair $(X,B)$ is Kawamata log terminal, 
and $D$ a Cartier divisor on $X$. Assume that $D$ is nef and 
that $D-(K_X+B)$ is nef and big. Then $H^0(X,D)\neq 0$.
\end{conj}

This conjecture is closely related to the minimal model program 
and plays an important role in the classification theory of Fano 
varieties. For a detailed introduction to this conjecture, 
we refer the reader to \cite{xie06}.

By the Kawamata-Viehweg vanishing theorem, 
we have $H^i(X,D)=0$ for any positive integer $i$. 
Thus $H^0(X,D)\neq 0$ is equivalent to $\chi(X,D)\neq 0$. 
Under the same assumptions as in Conjecture \ref{1}, 
the Kawamata-Shokurov non-vanishing theorem says 
that $H^0(X,mD)\neq 0$ for all $m\gg 0$. 
Thus the effective non-vanishing conjecture is an 
improvement of the non-vanishing theorem in some sense.

Note that $EN_1$ is trivial by the Riemann-Roch theorem, and 
that $EN_2$ was settled by Kawamata \cite[Theorem 3.1]{ka00} 
by virtue of the logarithmic semipositivity theorem. 
For $n\geq 3$, only a few results are known. 
For instance, $EN_n$ holds trivially for toric varieties 
\cite{mus}, $EN_3(X,0)$ holds for all canonical projective 
minimal threefolds $X$ \cite[Proposition 4.1]{ka00}, 
and $EN_3(X,0)$ also holds for almost all of canonical projective 
threefolds $X$ with $-K_X$ nef \cite[Corollary 4.5]{xie}.

In this note, we shall prove that, in the irregular case, 
the effective non-vanishing conjecture can be reduced to 
lower-dimensional cases by means of the Fourier-Mukai transform. 
As consequences, $EN_2$ is reproved after Kawamata, and $EN_n$ 
holds for all varieties of maximal Albanese dimension.

Throughout this note, we work over the complex number field $\C$. 
For the definition of Kawamata log terminal (KLT, for short) and 
the other notions, we refer the reader to \cite{kmm,km}.

For irregular varieties, the study of the Albanese map provides 
enough information to understand their birational structure. 
Therefore, through the Albanese map, we can utilize the Fourier-Mukai 
transform to give a reduction of the effective non-vanishing conjecture 
for irregular varieties. This idea was first used in \cite{ch}. 
First of all, we need the following lemma which follows easily from 
\cite[Theorem 2.2]{mu}.

\begin{lem}\label{2}
Let $A$ be an abelian variety, $\FF$ a coherent sheaf on $A$. 
Assume that $H^i(A,\FF\otimes P)=0$ for all $P\in\Pic^0(A)$ 
and all $i$. Then $\FF=0$.
\end{lem}

\begin{proof}
Let $\hat{A}$ be the dual abelian variety of $A$. The assumption 
implies that the Fourier-Mukai transform $\Phi(\FF)$ of $\FF$ is 
the zero sheaf on $\hat{A}$. Since the Fourier-Mukai transform 
$\Phi: D(A)\rightarrow D(\hat{A})$ induces an equivalence of 
derived categories \cite[Theorem 2.2]{mu}, we have $\FF=0$.
\end{proof}

\begin{thm}\label{3}
If $EN_k$ holds for any $k<n$, then $EN_n(X,B)$ holds for 
any $X$ with irregularity $q(X):=\dim H^1(X,\OO_X)>0$.
\end{thm}

\begin{proof}
By Kodaira's lemma, we may assume that $H=D-(K_X+B)$ is ample 
and $B$ is a $\Q$-divisor. Let $\pi: \widetilde{X}\rightarrow X$ 
be a resolution of $X$, and 
$\widetilde{\al}: \widetilde{X}\rightarrow A=\Alb(\widetilde{X})$ 
the Albanese morphism of $\widetilde{X}$. Since $(X,B)$ is KLT, 
$X$ has only rational singularities by \cite[Theorem 5.22]{km}, 
hence $q(\widetilde{X})=q(X)>0$. Since there are no rational curves 
on $A$, we have a non-trivial proper morphism $\al: X\rightarrow A$. 

Let $P\in\Pic^0(A)$, $P'=\al^*P$ and $\FF=\al_*\OO_X(D)$. 
By the Kawamata-Viehweg vanishing theorem, we have $H^i(X,D+P')=0$ 
for any $i>0$. By the relative Kawamata-Viehweg vanishing theorem 
\cite[Theorem 1-2-5]{kmm}, we have $R^i\al_*\OO_X(D+P')=0$ for any 
$i>0$. It follows from the Leray spectral sequence that 
$H^i(A,\FF\otimes P)=H^i(X,D+P')=0$ for any $i>0$. 
If $H^0(A,\FF)=0$, then $h^0(A,\FF\otimes P)=\chi(A,\FF\otimes P)=
\chi(A,\FF)=0$, i.e.\ $H^0(A,\FF\otimes P)=0$ for all $P\in\Pic^0(A)$. 
By Lemma \ref{2}, we have $\FF=0$.

Next we prove that $\FF\neq 0$, which implies 
$H^0(X,D)=H^0(A,\FF)\neq 0$. Let $a(X)=\dim\al(X)>0$. 
If $a(X)=n$, then $\al: X\rightarrow \al(X)$ is generically 
finite, and it is easy to see that $\FF\neq 0$. 
Assume that $a(X)<n$. 
Let $f:X\rightarrow Y$ be the Stein factorization of 
$\al$, $F$ a general fiber of $f$ and $\GG=f_*\OO_X(D)$. 
Then $F$ is a normal proper variety of dimension less than $n$. 
Note that $D|_F$ is nef Cartier, $(F,B|_F)$ is KLT and 
$D|_F-(K_F+B|_F)=H|_F$ is ample. By assumption, we have 
$\rk\GG=h^0(F,D|_F)\neq 0$, hence $\GG\neq 0$ as well as 
$\FF\neq 0$.
\end{proof}

\begin{cor}\label{4}
$EN_2$ holds, and $EN_3$ holds for any $X$ with $q(X)>0$.
\end{cor}

\begin{proof}
For $EN_2$, by the Riemann-Roch theorem, one has only to deal with 
the case where $X$ is a ruled surface over a smooth projective curve 
$C$ with $q(X)=g(C)\geq 2$. Since $EN_1$ holds, $EN_2$ also holds by 
Theorem \ref{3}. The second conclusion is obvious.
\end{proof}

\begin{cor}\label{5}
$EN_n(X,B)$ holds for any $X$ of maximal Albanese dimension.
\end{cor}

\begin{proof}
By assumption, $X$ is of maximal Albanese dimension, i.e.\ 
the Albanese morphism $\al:X\rightarrow A$ satisfies that 
$\dim\al(X)=\dim X=n$. So we can repeat the same argument 
as in Theorem \ref{3} to complete the proof by noting that 
$\al$ is generically finite.
\end{proof}

\begin{rem}\label{6}
Note that Corollary \ref{5} has already appeared in \cite{pp03} 
and \cite[Theorem 5.8]{pp05}. 
Note also that the assumption that $D$ is nef in Conjecture \ref{1} 
is not needed in the proof of Corollary \ref{5}, however 
\cite[Lemma 5.1]{pp05} proved that if $D-(K_X+B)$ is nef and big, 
then $D$ must be nef on the variety $X$ of maximal Albanese dimension. 
Furthermore, when the assumption that $D-(K_X+B)$ is nef and big 
in Conjecture \ref{1} is replaced with the weaker assumption that 
$D-(K_X+B)$ is either nef or of non-negative Iitaka dimension, 
$H^0(X,D)\neq 0$ also holds for any $X$ of maximal Albanese 
dimension \cite[Theorem 6.1]{pp06}. Finally, we should mention that 
Theorem \ref{3} and \cite[Theorem 5.8]{pp05} use a similar idea in proof.
\end{rem}

\begin{ack}
I am deeply indebted to Professors Yujiro Kawamata and 
Takao Fujita for valuable advices and warm encouragements. 
I would like to express my gratitude to Professors 
Giuseppe Pareschi and Mihnea Popa for pointing out 
some known results which are summarized in Remark \ref{6}.
I also thank the referee for useful suggestions and comments.
This work was partially supported by JSPS grant No.\ P05044 
and by the 21st Century COE Program.
\end{ack}

\textsc{Department of Mathematics, Tokyo Institute of Technology, 2-12-1 Oh-okayama, Meguro, Tokyo 152-8551, Japan}

\textit{Current address}: \textsc{Graduate School of Mathematical Sciences, 
University of Tokyo, Komaba, Meguro, Tokyo 153-8914, Japan}

\textit{E-mail address}: \texttt{xie\_qihong@hotmail.com}

\end{document}